\documentclass{amsart}
\usepackage{psfig}

\newtheorem{theorem}{Theorem}[section]
\newtheorem{lemma}[theorem]{Lemma}
\newtheorem{corollary}[theorem]{Corollary}

\theoremstyle{definition}

\theoremstyle{remark}

\numberwithin{equation}{section}

\newcommand{\abs}[1]{\lvert#1\rvert}

\DeclareSymbolFont{AMSb}{U}{msb}{m}{n}
\DeclareMathSymbol{\Z}{\mathalpha}{AMSb}{"5A}

\begin{document}
\newcommand{\beqs}{\begin{equation*}}
\newcommand{\eeqs}{\end{equation*}}
\newcommand{\beq}{\begin{equation}}
\newcommand{\eeq}{\end{equation}}
\newcommand\nutwid{\overset {\text{\lower 3pt\hbox{$\sim$}}}\nu}
\newcommand\Mtwid{\overset {\text{\lower 3pt\hbox{$\sim$}}}M}
\newcommand\ptwid{\overset {\text{\lower 3pt\hbox{$\sim$}}}p}
\newcommand\pitwid{\overset {\text{\lower 3pt\hbox{$\sim$}}}\pi}
\newcommand\bijone{\overset {1}\longrightarrow}                   
\newcommand\bijtwo{\overset {2}\longrightarrow}                   
\newcommand\pihat{\widehat{\pi}}
\newcommand\mymod[1]{(\mbox{mod}\ {#1})}
\newcommand\myto{\to}
\newcommand\srank{\mathrm{srank}}
\newcommand\pitc{\pi_{\mbox{$t$-core}}}
\newcommand\pitcb[1]{\pi_{\mbox{${#1}$-core}}}  
\newcommand\mpitcb[1]{\pi_{\mathrm{{#1}-core}}}  
\newcommand\lamseq{(\lambda_1, \lambda_2, \dots,\lambda_\nu)}
\newcommand\mylabel[1]{\label{#1}}
\newcommand\eqn[1]{(\ref{eq:#1})}
\newcommand\stc{{St-crank}}
\newcommand\mstc{\mbox{St-crank}}
\newcommand\mmstc{\mbox{\scriptsize\rm St-crank}}
\newcommand\tqr{{$2$-quotient-rank}}
\newcommand\mtqr{\mbox{$2$-quotient-rank}}
\newcommand\mmtqr{\mathrm{2-quotient-rank}}
\newcommand\fcc{{$5$-core-crank}}        
\newcommand\bgr{BG-rank} 
\newcommand\mbgr{\mbox{BG-rank}} 
\newcommand\nvec{(n_0, n_1, \dots, n_{t-1})}
\newcommand\parity{\mbox{par}}
\newcommand\pbar{\overline{p}}
\newcommand\Pbar{\overline{P}}
\newcommand\abar{\overline{a}}

\title[The Andrews-Stanley refinement of Ramanujan's congruence]{On the Andrews-Stanley Refinement of \\
Ramanujan's Partition Congruence \\
Modulo $5$ and Generalizations}

\author{Alexander Berkovich}
\address{Department of Mathematics, University of Florida, Gainesville,
Florida 32611-8105}
\email{alexb@math.ufl.edu}          

\author{Frank G. Garvan}
\address{Department of Mathematics, University of Florida, Gainesville,
Florida 32611-8105}
\email{frank@math.ufl.edu}          

\subjclass{Primary 11P81, 11P83; Secondary 05A17, 05A19}

\date{February 24, 2004}  


\keywords{partitions, $t$-cores, ranks, cranks, Stanley's statistic, Ramanujan's
congruences}

\begin{abstract}
In a recent study of sign-balanced, labelled posets Stanley \cite{Stan}, 
introduced a new integral partition statistic
\beqs
\mathrm{srank}(\pi) = {\mathcal O}(\pi) - {\mathcal O}(\pi'),
\eeqs
where ${\mathcal O}(\pi)$ denotes the number of odd parts of the partition $\pi$
and $\pi'$ is the conjugate of $\pi$. In \cite{Andrews1} Andrews proved
the following refinement of Ramanujan's partition congruence mod $5$:
\begin{align*}
p_0(5n+4) &\equiv p_2(5n+4) \equiv 0 \pmod{5},\\
p(n) &= p_0(n) + p_2(n),
\end{align*}
where $p_i(n)$ ($i=0,2$) denotes the number of partitions of $n$ with
$\mathrm{srank}\equiv i\pmod{4}$ and $p(n)$ is the number of unrestricted partitions
of $n$. Andrews asked for a partition statistic that would divide the
partitions enumerated by $p_i(5n+4)$ ($i=0,2$) into five equinumerous classes.

In this paper we discuss three such statistics: the \stc, the \tqr\ 
and the \fcc.
The first one, while new, is 
intimately related to the Andrews-Garvan \cite{AG} crank. 
The second one is in terms of the $2$-quotient of a partition.
The third one was introduced by Garvan, Kim and Stanton in \cite{GKS}.
We use it in our combinatorial proof of the Andrews refinement.
Remarkably, the Andrews result is a simple consequence of a stronger refinement
of Ramanujan's congruence mod $5$.  This more general refinement uses a new
partition statistic which we term the \bgr. We employ the \bgr\  to prove new
partition congruences modulo $5$.
Finally, we discuss some new formulas for partitions that are $5$-cores and
discuss an intriguing relation between $3$-cores and the Andrews-Garvan crank.
\end{abstract}

\maketitle

\section{Introduction} \label{sec:intro}
Let $p(n)$ be the number of unrestricted partitions of $n$. 
Ramanujan discovered and later proved that
\begin{align}
p(5n+4) &\equiv 0 \pmod{5},\mylabel{eq:ram5} \\ 
p(7n+5) &\equiv 0 \pmod{7},\mylabel{eq:ram7} \\ 
p(11n+6) &\equiv 0 \pmod{11}.\mylabel{eq:ram11} 
\end{align}
Dyson \cite{Dyson} was the first to consider combinatorial explanations
of these congruences. He defined the rank of a partition as 
{\it the largest part minus the number of parts} and made the
empirical observations that
\begin{align}
N(k,5,5n+4) &= \frac{p(5n+4)}{5},\quad 0 \le k \le 4, \mylabel{eq:dys5} \\ 
N(k,7,7n+5) &= \frac{p(7n+5)}{7},\quad 0 \le k \le 6, \mylabel{eq:dys7}  
\end{align}
where $N(k,m,n)$ denotes the number of partitions of $n$ with
rank congruent to $k$ modulo $m$.
Equation (\ref{eq:dys5}) means that the residue of
the rank mod ${5}$ divides the partitions of $5n+4$ 
into five equal classes. Similarly, 
(\ref{eq:dys7}) implies that the residue of
the rank mod ${7}$ divides the partitions of $7n+5$
into seven equal classes. Dyson's rank failed to explain (\ref{eq:ram11}),
and so Dyson conjectured the existence of a hypothetical statistic, called the
crank, that would explain the Ramanujan congruence mod ${11}$.
Identities (\ref{eq:dys5})-(\ref{eq:dys7}) were later proved by
Atkin and Swinnerton-Dyer \cite{ASD}. Andrews and Garvan \cite{AG} found a
crank for all three Ramanujan congruences (\ref{eq:ram5})-(\ref{eq:ram11}).
Their crank is defined as follows
\begin{equation}
\mbox{crank}(\pi) =
\begin{cases}
  \ell(\pi), &\mbox{if $\mu(\pi)=0$}, \\
  \nutwid(\pi) - \mu(\pi), &\mbox{if $\mu(\pi)>0$},
\end{cases}
\mylabel{eq:crankdef} 
\end{equation}
where $\ell(\pi)$ denotes the largest part of $\pi$,
$\mu(\pi)$ denotes the number of ones in $\pi$ and $\nutwid(\pi)$
denotes the number of parts of $\pi$ larger than $\mu(\pi)$.

Later, Garvan, Kim and Stanton \cite{GKS} found different cranks, which also
explained all three congruences (\ref{eq:ram5})-(\ref{eq:ram11}).
Their approach made essential use of $t$-cores of partitions 
and led to explicit bijections between various equinumerous classes.
In particular, they provided what amounts to a combinatorial proof
of the formula
\beq
\sum_{n\ge0} p(5n+4) q^n
= 5 \prod_{m\ge1} \frac{ (1-q^{5m})^5 }{ (1-q^m)^6 },
\mylabel{eq:rambest} 
\eeq
considered by Hardy to be an example of Ramanujan's best work.

The main results of \cite{AG} can be summarized as
\begin{align}
M(k,5,5n+4) &= \frac{p(5n+4)}{5},\quad 0 \le k \le 4, \mylabel{eq:ag5}\\ 
M(k,7,7n+5) &= \frac{p(7n+5)}{7},\quad 0 \le k \le 6, \mylabel{eq:ag7}\\ 
M(k,11,11n+6) &= \frac{p(11n+6)}{11},\quad 0 \le k \le 10, \mylabel{eq:ag11} 
\end{align}
and
\begin{align}
&1 + (x + x^{-1} - 1)q + \sum_{n>1} \sum_{m} \Mtwid(m,n) x^m q^n \nonumber\\
&\qquad\qquad= 
\prod_{n\ge1} \frac{ (1 - q^n) }{ (1 - x q^n) (1 - x^{-1} q^n) },
\mylabel{eq:crankgf} 
\end{align}
where $\Mtwid(m,n)$ denotes the number of partitions of $n$ with crank $m$
and $M(k,m,n)$ denotes the number of partitions of $n$ with
crank congruent to $k$ modulo $m$.

In \cite{G2} Garvan found a refinement of (\ref{eq:ram5})
\beq
M(k,2,5n+4) \equiv 0 \pmod{5}, \quad k=0,1 \mylabel{eq:gref5} 
\eeq
together with the combinatorial interpretation
\beq
M(2k+\alpha,10,5n+4) = \frac{M(\alpha,2,5n+4)}{5},\quad 0 \le k \le 4, 
\mylabel{eq:gref5a}\\ 
\eeq
with $\alpha=0,1$.

Recently, a very different refinement of (\ref{eq:ram5}) was given by Andrews
\cite{Andrews1}. Building on the work of Stanley \cite{Stan}, Andrews examined
partitions $\pi$ classified according to ${\mathcal O}(\pi)$ and
${\mathcal O}(\pi')$, where 
where ${\mathcal O}(\pi)$ denotes the number of odd parts of the partition $\pi$
and $\pi'$ is the conjugate of $\pi$. He used recursive relations to show that
\beq
G(z,y,q) := \sum_{n,r,s\ge0} S(n,r,s) q^n z^r y^s
=
\frac{ (-zyq;q^2)_\infty }{ (q^4;q^4)_\infty (z^2q^2;q^4)_\infty  
(y^2q^2;q^4)_\infty},
\mylabel{eq:rsgf} 
\eeq
where $S(n,r,s)$ denotes the number of partitions $\pi$ of $n$ with
${\mathcal O}(\pi)=r$,
${\mathcal O}(\pi')=s$, and
\begin{align}
(a;q)_\infty &= \lim_{n\to\infty} (a;q)_n, \mylabel{eq:aqdef} \\ 
(a;q)_n = (a)_n &=
\begin{cases}
1, &\mbox{if $n=0$},\\
\prod_{j=0}^{n-1}(1-aq^j), &\mbox{if $n>0$.}
\end{cases}
\mylabel{eq:aqndef} 
\end{align}
A direct combinatorial proof of (\ref{eq:rsgf}) was later given
by A. Sills \cite{Sills}, A.~J.~Yee \cite{Y} and C. Boulet \cite{Boulet}. 
Actually, C. Boulet proved a stronger version of (\ref{eq:rsgf})
with one extra parameter.
We define the Stanley rank of a partition
$\pi$ as
\beq
\mathrm{srank}(\pi) = {\mathcal O}(\pi) - {\mathcal O}(\pi').
\mylabel{eq:srankdef} 
\eeq
It is easy to see that
\beq
\mathrm{srank}(\pi) \equiv 0 \pmod{2},
\mylabel{eq:srankcong} 
\eeq
so that
\beq
p(n) = p_0(n) + p_2(n),
\mylabel{eq:p02} 
\eeq
where $p_i(n)$ ($i=0,2$) denotes the number of partitions of $n$ with
$\mathrm{srank}\equiv i\pmod{4}$.
We note that (\ref{eq:rsgf}) with $z=y^{-1}=\sqrt{-1}$ immediately implies
the Stanley formula \cite[p.8]{Stan}
\beq
\sum_{n\ge0} (p_0(n) - p_2(n)) q^n
= \frac{ (-q;q^2)_\infty }{ (q^4;q^4)_\infty (-q^2;q^4)_\infty^2}.
\mylabel{eq:p02prod} 
\eeq
Using (\ref{eq:ram5}), (\ref{eq:p02}) and (\ref{eq:p02prod}), Andrews
proved the following refinement of (\ref{eq:ram5})
\beq
p_0(5n+4) \equiv p_2(5n+4) \equiv 0 \pmod{5}.
\mylabel{eq:andrefine} 
\eeq
His proof of (\ref{eq:andrefine}) was analytic and so at the end of \cite{Andrews1}
he posed the problem of finding a partition statistic that would give a 
combinatorial interpretation of (\ref{eq:andrefine}).
The first goal of this paper is to provide such an interpretation. It turns out 
that there are several distinct integral partition statistics, whose residue
mod $5$ split the partitions enumerated by $p_i(5n+4)$ (with $i=0$, $2$) into five 
equal classes.  The first statistic,
which we call the \stc, is new. However, it is intimately related
to the Andrews-Garvan crank (\ref{eq:crankdef}).
The second statistic, which we call the \tqr, is also new. 
Unexpectedly, the third statistic is the \fcc,  introduced
by Garvan, Kim and Stanton \cite{GKS}.
This statistic not only provides the desired combinatorial
interpretation, but it also leads to a direct combinatorial proof of
(\ref{eq:andrefine}).

Our second goal here is to show that Andrews' result \eqn{andrefine}
is a straightforward corollary of the new refinement of \eqn{ram5}.
This stronger refinement uses a new partition statistic, which we term the \bgr.
Remarkably, the \bgr\ enables us to discover and prove new partition congruences
mod $5$.

The rest of this paper is organized as follows.
In Section 2 we define the \stc\  
and show that is indeed, a statistic asked for in \cite{Andrews1}. 
In Section 3 we give another combinatorial interpretation of \eqn{andrefine},
discuss a surprising relation between $3$-cores and the Andrews-Garvan crank,
and then briefly we review the development in \cite{GKS}. In Section 4 we 
establish a 
number of new formulas for partitions that are $5$-cores and outline a
combinatorial proof of \eqn{andrefine}. The hardest parts of this proof are 
dealt
with in Sections 5 and 6.
Finally, in Section 7 we introduce the \bgr\ and use it to prove
new partition congruences mod $5$.


\section{The \stc} \label{sec:stcrank}
We begin with some preliminaries about partitions and their conjugates.
A partition $\pi$ is a nonincreasing sequence
\beq
\pi = (\lambda_1, \lambda_2, \lambda_3, \dots)
\mylabel{eq:pidef} 
\eeq
of nonnegative integers (parts)
\beq
\lambda_1 \ge \lambda_2 \ge \lambda_3 \ge \cdots.
\mylabel{eq:lams} 
\eeq
The weight of $\pi$, denoted by $\abs{\pi}$ is the sum of parts
\beq
\abs{\pi} = \lambda_1 + \lambda_2 + \lambda_3 + \cdots.
\mylabel{eq:piweight} 
\eeq
If $\abs{\pi}=n$, then we say that $\pi$ is a partition of $n$. Often it is 
convenient to use another notation for $\pi$
\beq
\pi = (1^{f_1}, 2^{f_2}, 3^{f_3}, \dots),
\mylabel{eq:pidef2}
\eeq
which indicates the number of times each integer occurs as a part.
The number $f_i=f_i(\pi)$ is called the frequency of $i$ in $\pi$.
The conjugate of $\pi$ is the partition 
$\pi'=(\lambda_1', \lambda_2', \lambda_3', \dots)$ with
\begin{align}
\lambda_1' &= f_1 + f_2 + f_3 + f_4 + \cdots \nonumber\\
\lambda_2' &= f_2 + f_3 + f_4 + \cdots \mylabel{eq:piconj}\\
\lambda_3' &= f_3 + f_4 + \cdots \nonumber\\
	  &\vdots \nonumber
\end{align}
Next, we discuss two bijections. The first one relates $\pi$ and bipartitions
$(\pi_1,\pi_2)$, where $\pi_2$ is a partition with no repeated even parts.

\noindent
{\bf Bijection 1}  
\beqs
\pi \bijone (\pi_1,\pi_2),
\eeqs
where
\begin{align*}
\pi &= (1^{f_1}, 2^{f_2}, 3^{f_3}, \dots),\\
\pi_1 &= (1^{\lfloor{f_2/2}\rfloor}, 2^{\lfloor{f_4/2}\rfloor}, 
          3^{\lfloor{f_6/2}\rfloor}, \dots),\\
\pi_2 &= (1^{f_1}, 2^{\left\{f_2\right\}}, 3^{f_3}, 4^{\left\{f_2\right\}},
\dots),
\end{align*}
$\lfloor{x}\rfloor$ is the largest integer $\le x$, and
\beqs
\left\{x\right\} = x - 2 \lfloor{x/2}\rfloor.
\eeqs

\noindent
Indeed, remove from $\pi$ the maximum even number of even parts.
The resulting partition is $\pi_2$, The removed even parts can be organized
into a new partition
$(2^{2\lfloor{f_2/2}\rfloor}, 4^{2\lfloor{f_4/2}\rfloor},
6^{2\lfloor{f_6/2}\rfloor}, \dots),$
which can easily be mapped onto $\pi_1$.
Clearly, we have
\begin{align}
\abs{\pi} &= 4\abs{\pi_1} + \abs{\pi_2}, \mylabel{eq:bij1prop1}\\ 
\mathrm{srank}(\pi) &= \mathrm{srank}(\pi_2), \mylabel{eq:bij1prop2} 
\end{align}
so that
\begin{align}
\sum_{\pi} q^{\abs{\pi}} y^{\mathrm{srank}(\pi)}
&= \sum_{\pi_1} q^{4\abs{\pi_1}}
\sum_{\pi_2} q^{\abs{\pi_2}} y^{\mathrm{srank}(\pi_2)} \nonumber\\
&=\frac{1}{(q^4;q^4)_\infty} 
\sum_{\pi_2} q^{\abs{\pi_2}} y^{\mathrm{srank}(\pi_2)}.
\mylabel{eq:bij1gf} 
\end{align}
Comparing (\ref{eq:bij1gf}) and  (\ref{eq:rsgf}) with $zy=1$, we see that
\beq
\sum_{\pi_2} q^{\abs{\pi_2}} y^{\mathrm{srank}(\pi_2)}
= \frac{ (-q;q^2)_\infty }{ (y^2q^2;q^4)_\infty (q^2/y^2;q^4)_\infty },
\mylabel{eq:srankprodid} 
\eeq
where the sum is over all partitions with no repeated even parts.

To describe our second bijection we require a few definitions.
We say that $\pi_A$ is a partition of type A iff 
$\pi_A \bijone ((1),\pi_2)$. We say that 
$\pi_B = (\lambda_1, \lambda_2, \lambda_3, \dots)$ is
a partition of type B iff either $\abs{\pi_B}\ne4$,
$\lambda_1-\lambda_2\ge2$,
$\lambda_1'-\lambda_2'\ge2$,
$\lambda_1 -2$ and $\lambda_2$ are not identical even integers and
$\pi_B$ has no repeated even parts, or $\pi_B=(3,1)$.
Obviously, $\pi_B \bijone ((0),\pi_B)$.
Our second bijection relates partitions of type A and B.

\noindent
{\bf Bijection 2}  
\beqs
\pi_A \bijtwo \pi_B,
\eeqs
where
\begin{align*}
\pi_A &= (1^{f_1}, 2^{f_2}, 3^{f_3}, \dots, m^{f_m}),\\
\pi_B &=
\begin{cases}
(1^{f_1+2}, 2^{f_2-2}, 3^{f_3}, 4^{f_4}, \dots, (m-1)^{f_{m-1}}, m^{f_m-1},
(m+2)^1), &\mbox{if $m>2$}, \\
(1^{f_1+2},4^1), &\mbox{if $m=2$, $f_2=3$}, \\
(1^{f_1+1},3^1), &\mbox{if $m=2$, $f_2=2$},
\end{cases}
\end{align*}
$m\ge2$, $f_2=2$, $3$, and $f_{2i}=0$, $1$ for $i>1$.

\begin{figure}
\centerline{\psfig{figure=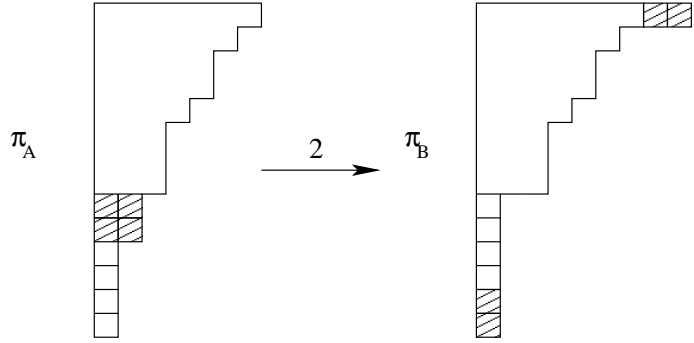}}
\caption{Graphical illustration of Bijection 2}
\label{fig1}
\end{figure}

Clearly, we have
\begin{align}
\abs{\pi_A} &= \abs{\pi_B}, \mylabel{eq:bij2prop1} \\ 
\mathrm{srank}(\pi_A) &= \mathrm{srank}(\pi_B). \mylabel{eq:bij2prop2} 
\end{align}
Next, we define a new partition statistic
\beq
\mstc(\pi) = \mathrm{crank}(\pi_1) + \frac{1}{2}\mathrm{srank}(\pi)
+ \Psi(\pi),
\mylabel{eq:stcrankdef} 
\eeq
where $\pi_1$ is determined by $\pi \bijone (\pi_1,\pi_2)$,
and the correction term $\Psi(\pi)=1$ if $\pi$ is of type B and zero, otherwise.
We note that
\beq
\mstc(\pi_A) = -1 + \frac{1}{2}\mathrm{srank}(\pi_A),
\mylabel{eq:stcrankprop1} 
\eeq
and
\beq
\mstc(\pi_B) = 1 + \frac{1}{2}\mathrm{srank}(\pi_B).
\mylabel{eq:stcrankprop2} 
\eeq

We give some examples.
Let $\pi=(1^2,2^4,3^2,4^1,5^1,6^2)$. Then 
$\pi \bijone ( (1^2,3^1)$, $(1^2,3^2,4^1,5^1))$ so that
$\mstc(\pi) = (1-2) + (5-1)/2 + 0 = 1$.
Next, we consider a partition of type B. Let $\pi_B = (1^2,3^1,5^1)$.
Then 
$\mstc(\pi_B)= 0 + (4-2)/2 + 1 = 2$.

Equipped with the definitions above, we can now prove the following lemma.

\begin{lemma}
\label{lemma1}
If
\beqs
g(x,y,q) := \sum_{\pi} q^{\abs{\pi}} x^{\mmstc(\pi)}
             y^{\mathrm{srank}(\pi)},
\eeqs
then $g(x,y,q)$ has the product representation
\beqs
g(x,y,q) = \frac{ (q^4;q^4)_\infty (-q;q^2)_\infty }
                {(q^4x,q^4/x, q^2y^2x, q^2/(y^2x); q^4)_\infty},
\eeqs
where
\beqs
(a_1,a_2,a_2,\dots ;q)_\infty = (a_1;q)_\infty (a_2;q)_\infty (a_3;q)_\infty
\cdots.
\eeqs
\end{lemma}
\begin{proof}
If $\pi$ is not of type B and $\pi \bijone (\pi_1,\pi_2)$, then
using (\ref{eq:bij1prop1})-(\ref{eq:bij1prop2}), (\ref{eq:stcrankdef})
we find that
\beq
q^{\abs{\pi}} x^{\mmstc(\pi)} y^{\mathrm{srank}(\pi)}
= q^{4\abs{\pi_1} + \abs{\pi_2}} x^{\mathrm{crank}(\pi_1)}
  (xy^2)^{\mathrm{srank}(\pi_2)/2}.
\mylabel{eq:qxy} 
\eeq
On the other hand, if $\pi=\pi_A$ and $\pi_A \bijtwo \pi_B$, then
\begin{align}
&q^{\abs{\pi_A}} x^{\mmstc(\pi_A)} y^{\mathrm{srank}(\pi_A)}
+
q^{\abs{\pi_B}} x^{\mmstc(\pi_B)} y^{\mathrm{srank}(\pi_B)}
\nonumber\\
&\quad=
q^{\abs{\pi_A}} (x+x^{-1}-1) (xy^2)^{\mathrm{srank}(\pi_A)/2}
+
q^{\abs{\pi_B}} x^0 (xy^2)^{\mathrm{srank}(\pi_B)/2}.
\mylabel{eq:qxy2} 
\end{align}
Here we have used (\ref{eq:bij1prop1})-(\ref{eq:bij1prop2})
and (\ref{eq:bij2prop1})-(\ref{eq:stcrankprop2}).

Equations (\ref{eq:qxy}), (\ref{eq:qxy2}) imply that
\beq
\sum_{\pi} q^{\abs{\pi}} x^{\mmstc(\pi)} y^{\mathrm{srank}(\pi)}
=
\sum_{\pi_1} q^{4\abs{\pi_1}} w(x,\pi_1)
\sum_{\pi_2} q^{\abs{\pi_2}} (xy^2)^{\mathrm{srank}(\pi_2)/2}
\mylabel{eq:stcrankgfid} 
\eeq
where
\beq
w(x,\pi_1) =
\begin{cases}
x+x^{-1}-1, &\mbox{if $\pi_1=(1)$,}\\
x^{\mathrm{crank}(\pi_1)}, &\mbox{otherwise.}
\end{cases}
\mylabel{eq:wdef} 
\eeq
We note that in the first sum on the right side of (\ref{eq:stcrankgfid})
the summation is over unrestricted partitions $\pi_1$, and in the
second sum the summation is over partitions $\pi_2$ with no repeated even parts.
Finally, recalling (\ref{eq:crankgf}) with $q\to q^4$ and (\ref{eq:srankprodid})
with $y^2\to xy^2$, we obtain
\beqs
\sum_{\pi} q^{\abs{\pi}} x^{\mmstc(\pi)} y^{\mathrm{srank}(\pi)}
=
\frac{ (q^4;q^4)_\infty }{ (xq^4,q^4/x;q^4)_\infty }
\cdot
\frac{ (-q;q^2)_\infty }{ (xy^2q^2,q^2/(xy^2);q^4)_\infty },
\eeqs
as desired.
\end{proof}

Next we show that
\begin{align}
&\mbox{the coefficient of $q^{5n+4}$ in $g(\xi,1,q)=0$},
\mylabel{eq:coeffz1}\\ 
&\mbox{the coefficient of $q^{5n+4}$ in $g(\xi,\sqrt{-1},q)=0$},
\mylabel{eq:coeffzi} 
\end{align}
where $\xi$ is a primitive fifth root of unity ($\xi^5=1$).
We use the method of \cite{G1}.
We need Jacobi's triple product identity 
\beq
\sum_{n=-\infty}^\infty z^n q^{n^2} = (q^2,-qz,-q/z;q^2)_\infty,
\mylabel{eq:jtp} 
\eeq
which implies that
\beq
(q^4;q^4)_\infty (-q;q^2)_\infty = (q^4,-q^3,-q;q^4)_\infty
= \sum_{n=-\infty}^\infty q^{2n^2+n} = \sum_{k\ge0}q^{T_k},
\mylabel{eq:jtpa} 
\eeq
and
\beq
(q^2\xi^2,q^2/\xi^2,q^2;q^2)_\infty = 
\frac{1}{1-\xi^2} \sum_{m\ge0} (-1)^m q^{2T_m} \xi^{-2m}(1-\xi^{4m+2}).
\mylabel{eq:jtpb} 
\eeq
Here $T_k = k(k+1)/2$. 
By Lemma \ref{lemma1} and equations
(\ref{eq:jtpa}) and (\ref{eq:jtpb}) we have
\begin{align}
&g(\xi,1,q) = \frac{ \sum_{k\ge0}q^{T_k} }
                   {(q^4\xi,q^4/\xi,q^2\xi,q^2/\xi;q^4)_\infty}
=
\frac{(q^2\xi^2,q^2/\xi^2,q^2;q^2)_\infty}
     { (q^{10};q^{10})_\infty}
\sum_{k\ge0}q^{T_k} \nonumber \\ 
&\quad
= \frac{1}{1-\xi^2} \frac{1}{(q^{10};q^{10})_\infty}
\sum_{k,m\ge0} (-1)^m q^{2T_m+T_k} \xi^{-2m}(1-\xi^{4m+2}).
\mylabel{eq:gz1id} 
\end{align}
Note that $2T_m + T_k \equiv 4 \pmod{5}$ iff $k\equiv m\equiv 2\pmod{5}$,
but then $1-\xi^{4m+2}=0$. This proves (\ref{eq:coeffz1}).
The proof of (\ref{eq:coeffzi}) is analogous.

Let $P_i(k,m,n)$ denote the number of partitions of $n$ with
$\mathrm{srank}\equiv i\pmod{4}$ and $\mstc\equiv k\pmod{m}$.
Clearly,
\begin{align}
\sum_{k=0}^4 \xi^k \sum_{n\ge0}P_0(k,5,n) q^n &=
\frac{ g(\xi,1,q) + g(\xi,\sqrt{-1},q)}{2},\mylabel{eq:P0gf} \\ 
\sum_{k=0}^4 \xi^k \sum_{n\ge0}P_2(k,5,n) q^n &=
\frac{ g(\xi,1,q) - g(\xi,\sqrt{-1},q)}{2},\mylabel{eq:P2gf}  
\end{align}
Combining (\ref{eq:coeffz1})-(\ref{eq:coeffzi}) and
(\ref{eq:P0gf})-(\ref{eq:P2gf}) we find that
\beq
\sum_{k=0}^4 \xi^k P_i(k,5,5n+4)=0, (\mbox{for $i=0,2$}),
\mylabel{eq:Piz} 
\eeq
which implies that
\beq
P_i(0,5,5n+4) = P_i(1,5,5n+4) = \cdots = P_i(4,5,5n+4).
\mylabel{eq:Pi04} 
\eeq
On the other hand
\beq
p_i(5n+4)=\sum_{k=0}^4 P_i(k,5,5n+4),
\mylabel{eq:pi54} 
\eeq
so that
\beq
P_i(k,5,5n+4) = \frac{1}{5} p_i(5n+4),
\mylabel{eq:mainresult} 
\eeq
for $i=0,2$ and $k=0,1,2,3,4$.
Thus, we have proved the main result of this section.
\begin{theorem}
\label{theorem1}
The residue of the partition statistic \stc\ mod $5$ divides the
partitions enumerated by $p_i(5n+4)$ with $i=0,2$ into
five equinumerous classes.
\end{theorem}
We illustrate this theorem in Table 1 below for the $30$
partitions of $9$. These partitions are organized into five classes
with six members each. In each class the first 4 members have
$\textrm{srank}\equiv0\pmod{4}$ and the remaining two members have
$\textrm{srank}\equiv2\pmod{4}$.

%

\begin{table}[ht]
\caption{}\label{table1}
\beqs
\renewcommand\arraystretch{1.1}
\begin{array}{|c|c|c|c|c|c|}
\hline
&\mstc\equiv0\mymod{5} &
1\mymod{5} &
2\mymod{5} &
3\mymod{5} &
4\mymod{5} \\
\hline
\mathrm{srank}\equiv0 &({3}^{3}) &
({1}^{5},{2}^{2}) &
({1}^{4},{2}^{1},{3}^{1}) &
({1}^{1},{2}^{4}) &
({1}^{9}) \\
\mymod{4} &({1}^{3},{2}^{1},{4}^{1}) &
({1}^{4},{5}^{1}) &
({1}^{3},{3}^{2}) &
({1}^{6},{3}^{1}) &
({1}^{2},{2}^{2},{3}^{1}) \\
& ({1}^{1},{3}^{1},{5}^{1}) &
({1}^{2},{2}^{1},{5}^{1}) &
({1}^{1},{4}^{2}) &
({1}^{1},{2}^{1},{6}^{1}) &
({2}^{3},{3}^{1}) \\
 & ({4}^{1},{5}^{1}) &
({9}^{1}) &
({2}^{2},{5}^{1}) &
({2}^{1},{7}^{1}) &
({1}^{2},{7}^{1}) \\                
\hline
\mathrm{srank}\equiv2 & ({1}^{3},{2}^{3}) &
({1}^{1},{2}^{1},{3}^{2}) &
({1}^{5},{4}^{1}) &
({1}^{7},{2}^{1}) &
({2}^{1},{3}^{1},{4}^{1}) \\
\mymod{4} &({1}^{3},{6}^{1}) &
({1}^{2},{3}^{1},{4}^{1}) &
({1}^{1},{8}^{1}) &
({1}^{1},{2}^{2},{4}^{1}) &
({3}^{1},{6}^{1}) \\
\hline
\end{array}
\eeqs
\end{table}

Finally, we note that the equation
\beq
\mathrm{srank}(\pi) = - \mathrm{srank}(\pi') \mylabel{eq:srconjprop} 
\eeq
implies that a partition $\pi$ is self-conjugate only if $\mathrm{srank}(\pi)=0$.
This means that the involution $\pi \longrightarrow \pi'$ has no fixed points
if $\mathrm{srank}(\pi)\equiv2\pmod{4}$. Hence, $2\mid p_2(5n+4)$ and
by (\ref{eq:andrefine}) we have the stronger congruence
\beqs
p_2(5n+4) \equiv 0 \pmod{10}.
\eeqs

\section{$t$-cores} \label{sec:tcores}

\subsection{Preliminaries} \label{subsec:pre}
In this section we recall some basic facts about $t$-cores and briefly
review the development in \cite{GKS}. A partition $\pi$ is called
a $t$-core, if it has no rim hooks of length $t$ \cite{JK}.
We let $a_t(n)$ denote the number of partitions of $n$ which are $t$-cores.
In what follows, $\pi_{\mbox{$t$-core}}$ denotes a $t$-core partition.
Given the diagram of a partition $\pi$ we label a cell
in the $i$-th row and $j$-th column by the least nonnegative integer congruent
to $j-i\pmod{t}$. The resulting diagram
is called a $t$-residue diagram \cite[p.84]{JK}. 

We also label cells in the infinite column $0$ and in the infinite row $0$ in the same 
way, and call the resulting diagram the extended $t$-residue diagram \cite{GKS}.
A region $r$ in the extended diagram is the set of cells $(i,j)$ satisfying
$t(r-1) \le j-i < tr$. A cell is called exposed if it is at the end of a row.
One can construct $t$ bi-infinite words $W_0$, $W_1$, \dots, $W_{t-1}$ of
two letters $N$ (not exposed) and $E$ (exposed) as follows:
\beqs
\mbox{The $j$-th element of $W_i$} =
\begin{cases}
N, &\mbox{if $i$ is not exposed in region $j$}, \\
E, &\mbox{if $i$ is exposed in region $j$}. 
\end{cases}
\eeqs

Let $P$ be the set of all partitions and $P_{\mbox{$t$-core}}$ be the
set of all $t$-cores. There is well-known bijection which goes back to
Littlewood \cite{L}.
$\phi_1 \,:\, P \rightarrow P_{\mbox{$t$-core}} \times P \times \cdots \times P$,
\begin{align}
\phi_1(\pi) &= (\pi_{\mbox{$t$-core}}, \vec{\pihat_t}), \mylabel{eq:phi1} \\ 
\vec{\pihat_t} &= (\pihat_0, \pihat_1, \pihat_2, \dots, \pihat_{t-1}),
\mylabel{eq:pihatvecdef} 
\end{align}
such that
\beq
\abs{\pi} = \abs{\pi_{\mbox{$t$-core}}} 
+ t \sum_{i=0}^{t-1}\abs{\pihat_i}.
\mylabel{eq:phi1a} 
\eeq
This bijection is described in more detail in \cite{JK}, \cite{GKS} and \cite{G3}.
The following identity is an immediate corollary of this bijection.
\beq
\frac{1}{(q)_\infty} = \sum_{n\ge0}p(n) q^n
= \frac{1}{(q^t;q^t)_\infty^t} \sum_{n\ge0} a_t(n) q^n.
\mylabel{eq:phi1cor} 
\eeq
It can be rewritten as
\beq
\sum_{n\ge0} a_t(n) q^n = \frac{(q^t;q^t)_\infty^t}{(q)_\infty}.
\mylabel{eq:tcoregfid} 
\eeq

There is another bijection $\phi_2$, introduced in \cite{GKS}.
It is for $t$-cores only.
$\phi_2\,:\, P_{\mbox{$t$-core}} \rightarrow \{\vec{n}
=(n_0, n_1, \dots, n_{t-1}) \,:\, n_i\in\Z, n_0+\cdots+n_{t-1}=0\}$, 
\beq
\phi_2(\pi_{\mbox{$t$-core}}) = \vec{n}=(n_0,n_1,n_2, \dots, n_{t-1}).
\mylabel{eq:phi2} 
\eeq
We call $\vec{n}$ an $n$-vector. It has the following properties.
\beq
\vec{n}\in\Z^t,\qquad \vec{n}\cdot \vec{1}_t = 0,
\mylabel{eq:phi2prop1} 
\eeq
and
\beq
\abs{\pi_{\mbox{$t$-core}}}= \frac{t}{2} \sum_{i=0}^{t-1} n_i^2
+ \sum_{i=0}^{t-1} i n_i,
\mylabel{eq:phi2prop2} 
\eeq
where the $t$-dimensional vector $\vec{1}_t$ has all components
equal to $1$. The generating function identity that corresponds to
this second bijection is

\beq
\sum_{n\ge0} a_t(n) q^n = 
\sum_{\substack{\vec{n}\in\Z^t \\ \vec{n}\cdot\vec{1}_t=0}}
q^{\frac{t}{2}||\vec{n}||^2 + \vec{b}_t\cdot\vec{n}}.
\mylabel{eq:tcoregfid2} 
\eeq

Here
\beq
||\vec{n}||^2=\sum_{i=0}^{t-1} n_i^2,\quad\mbox{and}\quad
\vec{b}_t=(0,1,2, \dots, t-1).
\mylabel{eq:nvecdef} 
\eeq
To construct the $n$-vector of $\pi_{\mbox{$t$-core}}$ in 
(\ref{eq:phi2}), we follow \cite{G3} and define
\beq
\vec{r}(\pi_{\mbox{$t$-core}}) = (r_0, r_1, r_2, \dots, r_{t-1}),
\mylabel{eq:rvec} 
\eeq
where for $0\le i\le t-1$, $r_i(\pi_{\mbox{$t$-core}})$ denotes the number
of cells labelled $i\pmod{t}$ in the $t$-residue diagram of 
$\pi_{\mbox{$t$-core}}$.  Then (\ref{eq:phi2}) can be given explicitly
as
\beq
\phi_2(\pi_{\mbox{$t$-core}}) = \vec{n} = 
(r_0-r_1, r_1-r_2, r_2-r_3, \dots, r_{t-1}-r_0).
\mylabel{eq:rvec2} 
\eeq

It was shown in \cite{GKS} that a partition is a $t$-core with $n$-vector
$\nvec$ iff for all $i=0$, \dots, $t-1$ the bi-infinite word $W_i$ is of
the form
\beqs
\begin{matrix}
\mbox{Region}: &\cdots\cdots\cdots & n_{i}-1 & n_i & n_{i}+1 & n_{i}+2 & \cdots\cdots\cdots \\
W_i: &\cdots\cdots\cdots & E       & E   & N       & N       & \cdots\cdots\cdots 
\end{matrix}
\eeqs
We note that $\tfrac{t}{2}||\vec{n}||^2$ is a multiple of $t$ since
$\vec{n}\cdot\vec{1}_t=0$. Hence by (\ref{eq:phi1cor}) and (\ref{eq:tcoregfid2}) 
we have
\beq
\sum_{n\ge0} a_t(tn+\delta) q^{tn+\delta}
=
\sum_{\substack{\vec{n}\in\Z^t,\ \vec{n}\cdot\vec{1}_t=0 \\
\vec{n}\cdot\vec{b}_t\equiv \delta\pmod{t}}}
q^{\frac{t}{2}||\vec{n}||^2 + \vec{b}_t\cdot\vec{n}},
\mylabel{eq:tcoresift} 
\eeq
and
\beq
\sum_{n\ge0} p(tn+\delta) q^{n}
=
\frac{1}{(q)_\infty^t}
\sum_{n\ge0} a_t(tn+\delta) q^{n},
\mylabel{eq:psift} 
\eeq
where
$\delta=0,1,2$, \dots, $t-1$.

\subsection{The \tqr} \label{subsec:2quotrank}
Having collected the necessary background on $t$-cores, we are now in a position
to provide another combinatorial intrepretation of \eqn{andrefine}. To this
end we introduce the following new partition statistic
\beq
\mtqr(\pi) = \nu(\pihat_0) - \nu(\pihat_1),
\mylabel{eq:tqrdef}
\eeq
where $\pihat_0$ and $\pihat_1$ are determined by
\beq
\phi_1(\pi) = (\pitcb{2}, (\pihat_0,\pihat_1)),
\mylabel{eq:pi01}
\eeq
where $\nu(\pi)$ denotes the number of parts of $\pi$.
Our main result here is
\begin{theorem}
\label{theorem2}
The residue of the \tqr\ mod $5$ divides the
partitions enumerated by $p_i(5n+4)$ with $i=0,2$ into
five equal classes.
\end{theorem}
\begin{proof}
We start by recalling Proposition 3.1(d) in \cite{Stan}:
\beq
\srank(\pi) \equiv |\pi| - |\pitcb{2}| \pmod{4}.
\mylabel{eq:sr2core} 
\eeq
Using \eqn{phi1a} with $t=2$, we obtain from \eqn{sr2core}
\beq
\srank(\pi) \equiv 2(|\pihat_0| + |\pihat_1|) \pmod{4}.
\mylabel{eq:sr2q}
\eeq
Next, we define the generating function
\beq
G_2(x,y,q) := 
\sum_{\pi} q^{\abs{\pi}} x^{\mmtqr(\pi)} y^{\mathrm{srank}(\pi)}
\mylabel{eq:G2def}
\eeq
It is possible to find a product representation for $G_2(x,\omega,q)$
when $\omega^4=1$, namely         
\begin{align}
&G_2(x,\omega,q)
=
\sum_{\substack{\mpitcb{2}\\ \pihat_0, \pihat_1}}
q^{|\mpitcb{2}| +2(|\pihat_0| + |\pihat_1|)} x^{\nu(\pihat_0) - \nu(\pihat_1)}
\omega^{2(|\pihat_0| + |\pihat_1|)} \mylabel{eq:G2prod} \\
& = \frac{\sum_{k\ge0} q^{T_k}}{(xq^2\omega^2,q^2\omega^2/x;q^2\omega^2)_\infty}
= \frac{ (q^4;q^4)_\infty (-q;q^2)_\infty}
{(xq^2\omega^2,q^2\omega^2/x;q^2\omega^2)_\infty}.
\nonumber
\end{align}
Here we have used \eqn{jtpa} along with the fact that a partition is a $2$-core
if and only if it is a staircase \cite{JK}.
Hence using
\beq
(x^{\pm1} q^2 \omega^2 ; q^2\omega^2)_\infty =
(x^{\pm1} q^2 \omega^2, x^{\pm1} q^4 ; q^4)_\infty,     
\mylabel{eq:xwprod}
\eeq
we find that
\beq
G_2(x,\sqrt{\pm1},q) =
\frac{(q^4;q^4)_\infty (-q;q^2)_\infty}
{(\pm xq^2,xq^4,\pm q^2/x,q^4/x;q^4)_\infty}.
\mylabel{eq:G2prod2}
\eeq
Comparing this with the product in Lemma \ref{lemma1}, we see that
\beq
G_2(x,\sqrt{\pm1},q) = g(x,\sqrt{\pm1},q).
\mylabel{eq:G2gid}
\eeq
This means that
\beq
\ptwid_i(m,n) = p_i^*(m,n),\qquad i=0,2.
\mylabel{eq:stctqr}
\eeq
Here $\ptwid_i(m,n)$ (resp. $p_i^*(m,n)$) denotes the number of
partitions of $n$ with $\srank\equiv i\pmod{4}$
and \tqr $=m$ (resp. \stc $=m$).
Theorem \ref{theorem2} follows easily from \eqn{stctqr} and Theorem \ref{theorem1}.
\end{proof}
We remark that the \stc\ and the \tqr\ are distinct statistics.
For example, $\mstc((5,4,1))=0$ but $\mtqr((5,4,1))=1$.
It would be interesting to find a direct combinatorial proof of \eqn{stctqr}.

\subsection{$3$-cores and the crank} \label{subsec:3cores}       
It is natural to attempt to extend the construction of the previous section
to ($3$-core, $3$-quotient). To this end we define
\beq
G_3(x,q) :=
\sum_{\pi} q^{\abs{\pi}} (x^{3n_1}+x^{3n_2+1} + x^{-3n_2-1})
x^{3(\nu(\pihat_1) - \nu(\pihat_2))},
\mylabel{eq:G3def}
\eeq
where
\begin{align}
\phi_1(\pi) &= (\pitcb{3}, (\pihat_0,\pihat_1,\pihat_2)), \mylabel{eq:3ca}\\
\phi_2(\pitcb{3}) &= (-n_1-n_2,n_1,n_2).\mylabel{eq:3cb}
\end{align}
Hence using \eqn{phi1a}, \eqn{phi2prop2} with $t=3$, we find that
\beq
|\pi| = Q_3(n_1,n_2) + 3(|\pihat_0| + |\pihat_1| + |\pihat_2|),
\mylabel{eq:pin3}
\eeq
where
\beq
Q_3(n_1,n_2) = 3(n_1^2 + n_1 n_2 + n_2^2) + n_1 + 2 n_2.
\mylabel{eq:Q3def}
\eeq
Clearly,
\begin{align}
&G_3(x,q) \mylabel{eq:G3id}\\
&=\frac{\displaystyle 
\sum_{n_1,n_2} q^{Q_3(n_1,n_2)} x^{3n_1} 
+
\sum_{n_1,n_2} q^{Q_3(n_1,n_2)} x^{3n_2+1} 
+
\sum_{n_1,n_2} q^{Q_3(n_1,n_2)} x^{-1-3n_2} 
}
{
(q^3,x^3q^3,q^3/x^3;q^3)_\infty
}.
\nonumber 
\end{align}
Next, we change summation variables in the first, second, and third sums respectively
as
\begin{align*}
 3n_1 &=n-m, \qquad 3n_2 = n + 2m,\\
 3n_1+1 &=-2n-m, \qquad 3n_2+1 = n - m,\\
 3n_1+1 &=-n-2m, \qquad 3n_2+1 = m - n,   
\end{align*}
respectively. In this way, we get
\beq
G_3(x,q) =
\frac{
\displaystyle
\sum_{n,m} q^{n^2 + m n + m^2} x^{n-m}
}
{
(q^3,x^3q^3,q^3/x^3;q^3)_\infty
}.
\mylabel{eq:G3id2}
\eeq
Remarkably, the numerator on the right side of \eqn{G3id2}
has a product representation. By equation (1.23) in \cite{HGB} 
we have
\beq
\sum_{n,m} q^{n^2 + m n + m^2} x^{n-m}
= (x + 1 + 1/x) (q;q)_\infty (q^3;q^3)_\infty
\frac{
(x^3q^3;q^3)_\infty (q^3/x^3;q^3)_\infty
}
{
(xq;q)_\infty (q/x;q)_\infty
}.
\mylabel{eq:HGBid}
\eeq
We have
\beq
G_3(x,q) =  (x + 1 + 1/x) 
\frac{
(q)_\infty
}
{
(xq;q)_\infty (q/x;q)_\infty
}.
\mylabel{eq:G3id3}
\eeq
Recalling \eqn{crankgf}, we see that \eqn{G3id3} gives a surprising
relation between $3$-cores and the Andrews-Garvan crank.
This certainly warrants further investigation.

\subsection{$5$-cores} \label{subsec:5cores}       
We now assume $t=5$. For the case $\delta=4$ the right side of 
(\ref{eq:tcoresift}) can be simplified using the the following change
of variables.
\begin{align}
n_0 &=\alpha_0 + \alpha_4, \nonumber \\
n_1 &=-\alpha_0 + \alpha_1 + \alpha_4, \nonumber \\
n_2 &=-\alpha_1 + \alpha_2, \mylabel{eq:ntoa} \\ 
n_3 &=-\alpha_2 + \alpha_3 - \alpha_4, \nonumber \\
n_4 &=-\alpha_3 - \alpha_4, \nonumber 
\end{align}
We find $\vec{n}$ is an $n$-vector satisfying 
$\vec{n}\cdot\vec{b}_5\equiv4\pmod{5}$ if and only if
\beq
\vec{\alpha} = (\alpha_0, \alpha_1, \alpha_2, \alpha_3, \alpha_4)\in\Z^5
\mylabel{eq:avec1} 
\eeq
and
\beq
\alpha_0 + \alpha_1 + \alpha_2 + \alpha_3 + \alpha_4 = 1.
\mylabel{eq:avec2} 
\eeq
We call $\vec{\alpha}$ and $\alpha$-vector.
Hence, by (\ref{eq:tcoresift}) and (\ref{eq:psift}) we have
\beq
\sum_{n\ge0} a_5(5n+4) q^{n+1}
=
\sum_{\substack{\vec{\alpha}\cdot\vec{1}_5=1 \\
\vec{\alpha}\in\Z^5}}
q^{Q(\vec{\alpha})},
\mylabel{eq:5coresift} 
\eeq
and
\beq
\sum_{n\ge0} p(5n+4) q^{n+1}
=
\frac{1}{(q)_\infty^5}
\sum_{\substack{\vec{\alpha}\cdot\vec{1}_5=1 \\
\vec{\alpha}\in\Z^5}}
q^{Q(\vec{\alpha})},
\mylabel{eq:psift5} 
\eeq
where
\beq
Q(\vec{\alpha}) = ||\vec{\alpha}||^2 - 
                (\alpha_0\alpha_1 + \alpha_1\alpha_2 + \cdots + \alpha_4\alpha_0),
\mylabel{eq:Qadef} 
\eeq
If $\abs{\pi}\equiv4\pmod{5}$ and $t=5$, we can combine 
bijections $\phi_1$ and $\phi_2$ into a single bijection
\beq
\Phi(\pi) = (\vec{\alpha},\vec{\pihat_5}),
\mylabel{eq:Phi} 
\eeq
such that
\beq
\abs{\pi} = 5Q(\vec{\alpha}) - 1 + 5 \sum_{i=0}^4 \abs{\pihat_i}.
\mylabel{eq:Phi1} 
\eeq
Next, following \cite{GKS} we define the $5$-core crank of $\pi$ when
$\abs{\pi}\equiv4\pmod{5}$ as
\beq
c_5(\pi)
=1 + \sum_{i=0}^4 i \alpha_i
\equiv 2(1 + n_0 - n_1 - n_2 + n_3)
\equiv 2 + \sum_{i=-2}^2 i r_{2-i} \pmod{5},
\mylabel{eq:5ccrank} 
\eeq
where $\alpha$ is determined by (\ref{eq:Phi}).

It is easy to check that $Q(\vec{\alpha})$ in ({\ref{eq:Phi1}) remains invariant 
under the following cyclic permutation
\beq
\widehat{C}_1(\vec{\alpha}) = (\alpha_4,\alpha_0,\alpha_1,\alpha_2,\alpha_3),
\mylabel{eq:cycperm} 
\eeq
while $c_5(\pi)$ increases by $1\pmod{5}$ under the map
\beq
\widehat{O}(\pi) = \Phi^{-1}(\widehat{C}_1(\vec{\alpha}), \vec{\pihat}_5).
\mylabel{eq:perm} 
\eeq
In other words, if $\abs{\pi}\equiv4\pmod{5}$, then
\beq
\abs{\pi} = \abs{\widehat{O}(\pi)},
\mylabel{eq:permprop1} 
\eeq
and
\beq
c_5(\pi)+1\equiv c_5(\widehat{O}(\pi))\pmod{5}.
\mylabel{eq:permprop2} 
\eeq
This suggests that all partitions of $5n+4$ can be organized into orbits.
Each orbit consists of five distinct members:
\beq
\pi,\ \widehat{O}(\pi),\ \widehat{O}^2(\pi),\ \widehat{O}^3(\pi),\ 
\widehat{O}^4(\pi),
\mylabel{eq:orb} 
\eeq
and each element of the orbit has a distinct $5$-core crank (mod $5$).
Clearly, the total number of such orbits is $\tfrac{1}{5} p(5n+4)$,
and so $p(5n+4)\equiv0\pmod{5}$. This summarizes the combinatorial
proof of ({\ref{eq:ram5}) given in \cite{GKS}.
If we apply the map $\widehat{O}$ (\ref{eq:perm}) to the partitions of 
$5n+4$ that are $5$-cores, we find that
\beq
a_5^0(5n+4)= a_5^1(5n+4)= \cdots = a_5^4(5n+4),
\mylabel{eq:5corerels} 
\eeq
where, for $0\le j\le4$,  $a_5^j(n)$ denotes the number of partitions of $n$
that are $5$-cores with $5$-core crank congruent to $j$ modulo $5$.
Hence,
\beq
a_5^j(5n+4)= \frac{1}{5} a_5(5n+4),\quad j=0,1,\dots,4,
\mylabel{eq:5corerels2} 
\eeq
which proves that
\beq
a_5(5n+4)\equiv0\pmod{5}.
\mylabel{eq:5corecong} 
\eeq
Actually, more is true. We have 
\beq
a_5(5n+4) = 5a_5(n).
\mylabel{eq:5corerel} 
\eeq
We sketch the combinatorial proof of (\ref{eq:5corerel}) given in \cite{GKS}.
See also \cite{G3}. The map
$\theta:\,P_{\mbox{$5$-core}}(n) \longrightarrow P_{\mbox{$5$-core}}^0(5n+4)$,
defined in terms of $n$-vectors as
\begin{align}
\vec{n} \mapsto \vec{n}'
&=(n_1+2n_2+2n_4+1,
-n_1-n_2+n_3+n_4+1,
2n_1+n_2+2n_3, \nonumber \\
&\qquad -2n_2-2n_3-n_4-1,
-2n_1-n_3-2n_4-1),
\mylabel{eq:theta} 
\end{align}
is a bijection. Here $P_{\mbox{$5$-core}}(n)$ is the set of all $5$-cores of $n$,
and $P_{\mbox{$5$-core}}^0(n)$ is set of all $5$-cores of $n$ with $5$-core
crank congruent to zero modulo $5$. Since $\theta$ is a bijection, we have
\beq
a_5(n)=a_5^0(5n+4).
\mylabel{eq:5corerel2} 
\eeq
The proof of (\ref{eq:5corerel}) easily follows from 
(\ref{eq:5corerels2}) and (\ref{eq:5corerel2}).
Finally, we remark that Ramanujan's result (\ref{eq:rambest}) is
a straightforward consequence of (\ref{eq:psift}) with$(t,\delta)=(5,4)$,
(\ref{eq:5corerel}), and (\ref{eq:tcoregfid}) with $t=5$.

\section{Refinement of Ramanujan's mod $5$ congruence, the srank and
the $5$-core crank} \label{sec:refinement}

In the previous section we discussed the combinatorial proof in \cite{GKS} of
Ramanujan's congruence (\ref{eq:ram5}) using the the $5$-core crank 
(\ref{eq:5ccrank}). It is somewhat unexpected that the $5$-core crank can
be employed to prove the refinement (\ref{eq:andrefine}) as well.

In fact, we were amazed to discover the following elegant formulas
\begin{align}
\mathrm{srank}(\pi_{\mbox{$5$-core}}) &\equiv \sum_{i=0}^4 (n_i + i)^3 \pmod{4}, 
\mylabel{eq:elegant1} \\ 
\mathrm{srank}(\pi) &\equiv \mathrm{srank}(\pi_{\mbox{$5$-core}}) 
    + \sum_{i=0}^4 \mathrm{srank}(\pihat_i) \mylabel{eq:elegant2} \\ 
&\quad + 2 \sum_{i=0}^4 \abs{\pihat_i}(n_i + i) \pmod{4},
\nonumber
\end{align}
where $\pi_{\mbox{$5$-core}}$, 
$\vec{\pihat}=(\pihat_0, \pihat_1, \pihat_2, \pihat_3, \pihat_4)$ are determined 
by (\ref{eq:phi1})
with $t=5$, and
\beqs
\vec{n}=(n_0,n_1,\dots,n_4)=\phi_2(\pi_{\mbox{$5$-core}}).
\eeqs
In spite of their simple appearance, the above formulas are far from obvious.
In Sections 4 and 5 we prove generalizations of 
(\ref{eq:elegant1})-(\ref{eq:elegant2}).
Here we restrict our attention to some implications of 
(\ref{eq:elegant1})-(\ref{eq:elegant2}).

First, we note that if $\abs{\pi_{\mbox{$5$-core}}}\equiv4\pmod{5}$, then
(\ref{eq:elegant1}) can be written in terms of an 
$\alpha$-vector (\ref{eq:ntoa}) as
\beq
\mathrm{srank}(\pi_{\mbox{$5$-core}})
\equiv \alpha_0\alpha_1(\alpha_0-\alpha_1) + \alpha_1\alpha_2(\alpha_1-\alpha_2)
+ \cdots +
\alpha_4\alpha_0(\alpha_4-\alpha_0) \pmod{4}.
\mylabel{eq:sravec} 
\eeq
Similarly, if $\abs{\pi}\equiv4\pmod{5}$, then
\begin{align}
\mathrm{srank}(\pi) &\equiv
\alpha_0\alpha_1(\alpha_0-\alpha_1) + \cdots +
\alpha_4\alpha_0(\alpha_4-\alpha_0)
+ \sum_{i=0}^4 \mathrm{srank}(\pihat_i) \nonumber \\
& \quad + 2\{ (\alpha_0+\alpha_4)\abs{\pihat_0} + (\alpha_2+\alpha_3)\abs{\pihat_1}
+ (\alpha_1+\alpha_2)\abs{\pihat_2} \mylabel{eq:sravec2} \\ 
& \quad + (\alpha_0+\alpha_1)\abs{\pihat_3} + (\alpha_3+\alpha_4)\abs{\pihat_4}\}
\pmod{4}.
\nonumber
\end{align}
Remarkably, (\ref{eq:sravec}) suggests that $\mathrm{srank}(\pi_{\mbox{$5$-core}})$
with $\abs{\pi_{\mbox{$5$-core}}}\equiv4\pmod{5}$ remains invariant
mod $4$ under the cyclic permutation (\ref{eq:cycperm}), and we have
the following refinement of (\ref{eq:5corerels2}):
\beq
a_{5,i}^j(5n+4)= \frac{1}{5} a_{5,i}(5n+4), 
\mylabel{eq:5corerelrefine} 
\eeq
where $j=0$,\dots,$4$ and $i=0$, $2$. Here $a_{5,i}(n)$ denotes the number of 
$5$-cores of $n$ with $\mathrm{srank}\equiv i\pmod{4}$, and
$a_{5,i}^j(n)$ denotes the number of
$5$-cores of $n$ with $\mathrm{srank}\equiv i\pmod{4}$ and $5$-core crank 
$\equiv j\pmod{5}$. Moreover, it is not difficult to verify that the map $\theta$,
given by (\ref{eq:theta}), preserves the srank mod $4$. Indeed, recalling
that $n_0+n_1+n_2+n_3+n_4=0$ we find after some simplication that
\begin{align}
\sum_{i=0}^4 ( (n_i+i)^3 - (n_i'+i)^3) 
&\equiv 2( n_0n_2(n_0+n_2) + n_1n_3(n_1+n_3)+ n_2n_3(n_2+n_3) \nonumber \\
& \quad + n_1(n_1+1) + n_2(n_2+1) + n_3(n_3+1)) \mylabel{eq:invarmod4} \\ 
& \equiv 0 \pmod{4},
\nonumber
\end{align}
where $\vec{n}'$ is defined in (\ref{eq:theta}). Hence, (\ref{eq:5corerel2}) and 
(\ref{eq:5corerel})
can be refined as
\beq
a_{5,i}(n)=a_{5,i}^0(5n+4), \qquad (i=0,2), \mylabel{eq:5corerelrefine2} 
\eeq
and
\beq
a_{5,i}(5n+4)=5 a_{5,i}(n), \qquad (i=0,2), \mylabel{eq:5corerelrefine3} 
\eeq
respectively.

It is less trivial to prove the $5$-core crank analogue of 
Theorem \ref{theorem1}. Namely,
\begin{theorem}
\label{theorem3}
The residue of the $5$-core crank  mod $5$ divides the
partitions enumerated by $p_i(5n+4)$ with $i=0,2$ into
five equal classes.
\end{theorem}
\begin{proof}
We sketch a proof using (\ref{eq:sravec2}) and (\ref{eq:5corerelrefine}).
We define the cyclic shift operator $\widehat{C}_2$ by
\beq
\widehat{C}_2(\vec{\pihat}_5) = (\pihat_4,\pihat_2,\pihat_3,\pihat_0,\pihat_1),
\mylabel{eq:C2def} 
\eeq
Next, we use (\ref{eq:C2def}) to modify (\ref{eq:perm}) as
\beq
\widehat{O}_s(\pi) = \Phi^{-1}(\widehat{C}_1(\vec{\alpha}), 
                               \widehat{C}_2(\vec{\pihat}_5)),
\mylabel{eq:newperm} 
\eeq
where
$\Phi(\pi)=(\vec{\alpha},\vec{\pihat}_5)$. Fix $i=0,2$.
By (\ref{eq:sravec2}) we see that $\widehat{O}_s$ preserves the srank mod $4$,
and we may assemble all partitions of $5n+4$ with
$\mathrm{srank}\equiv i\pmod{4}$ into orbits:
\beqs
\pi,\ \widehat{O}_s(\pi),\ \widehat{O}_s^2(\pi),\ \widehat{O}_s^3(\pi),\
\widehat{O}_s^4(\pi),
\eeqs
where $\pi$ is some partition of $5n+4$ with $\mathrm{srank}(\pi)\equiv i\pmod{4}$.
As before, each orbit contains exactly five members and the $5$-core
crank increases by $1$ mod $5$ along the orbit. The number of these
orbits is $\tfrac{1}{5} p_i(5n+4)$, consequently $p_i(5n+4)\equiv0\pmod{5}$ and 
the result follows.
\end{proof}

Theorem \ref{theorem3} is illustrated below in Table \ref{table2},
which contains all $30$ partitions of $9$, organized into $6$ orbits.
Each row in this table represents an orbit, and the first row lists
all partitions of $9$ that are $5$-cores. In the table we have also included
the image of each partition under the bijection $\phi_1$. Instead of
giving the full $5$-quotient we have used a short-hand notation.
Terms in the table have the form $\pi\to(\pi_{\mbox{$5$-core}},k)$
where $k$ indicates that $\vec{\pihat}_5=(\pihat_0, \pihat_1, \pihat_2,
 \pihat_3, \pihat_{4})$, where $\pihat_i=(1)$ if $i=k$, and $(0)$ otherwise.


{\footnotesize
\begin{table}[ht]
\caption{}\label{table2}
\beqs
\renewcommand\arraystretch{1.2}
\begin{array}{|c|c|c|c|c|c|}
\hline
&c_5\equiv0\mymod{5} &   
1\mymod{5} &
2\mymod{5} &
3\mymod{5} &
4\mymod{5} \\
\hline
\mathrm{srank}\equiv0& 
({1}^{4},{5}^{1}) & 
({1}^{3},{3}^{2}) & 
({1}^{4},{2}^{1},{3}^{1}) & 
({1}^{1},{2}^{1},{6}^{1}) & 
({2}^{2},{5}^{1}) \\ 
\mymod{4} & &&&& \\
\cline{2-6}
& ({1}^{5},{2}^{2})\to & 
({2}^{3},{3}^{1})\to & 
({1}^{2},{7}^{1})\to & 
({4}^{1},{5}^{1})\to & 
({1}^{3},{2}^{1},{4}^{1})\to \\ 
 & 
(({2}^{2}),3) & 
(({1}^{4}),2) & 
(({1}^{2},{2}^{1}),1) & 
(({1}^{1},{3}^{1}),4) & 
(({4}^{1}),0) \\ 
\cline{2-6}
 & 
({3}^{3})\myto & 
({1}^{9})\myto & 
({1}^{1},{3}^{1},{5}^{1})\myto & 
({1}^{2},{2}^{2},{3}^{1})\myto & 
({9}^{1})\myto \\ 
 & 
(({2}^{2}),2) & 
(({1}^{4}),1) & 
(({1}^{2},{2}^{1}),4) & 
(({1}^{1},{3}^{1}),0) & 
(({4}^{1}),3) \\ 
\cline{2-6}
 & 
({2}^{1},{7}^{1})\myto & 
({1}^{2},{2}^{1},{5}^{1})\myto & 
({1}^{1},{2}^{4})\myto & 
({1}^{6},{3}^{1})\myto & 
({1}^{1},{4}^{2})\myto \\ 
 & 
(({2}^{2}),1) & 
(({1}^{4}),4) & 
(({1}^{2},{2}^{1}),0) & 
(({1}^{1},{3}^{1}),3) & 
(({4}^{1}),2) \\ 
\hline
\mathrm{srank}\equiv2
 & 
({1}^{3},{2}^{3})\myto & 
({1}^{3},{6}^{1})\myto & 
({2}^{1},{3}^{1},{4}^{1})\myto & 
({1}^{1},{8}^{1})\myto & 
({1}^{2},{3}^{1},{4}^{1})\myto \\ 
\mymod{4} & 
(({2}^{2}),4) & 
(({1}^{4}),0) & 
(({1}^{2},{2}^{1}),3) & 
(({1}^{1},{3}^{1}),2) & 
(({4}^{1}),1) \\ 
\cline{2-6}
 & 
({3}^{1},{6}^{1})\myto & 
({1}^{1},{2}^{2},{4}^{1})\myto & 
({1}^{7},{2}^{1})\myto & 
({1}^{1},{2}^{1},{3}^{2})\myto & 
({1}^{5},{4}^{1})\myto \\ 
 & 
(({2}^{2}),0) & 
(({1}^{4}),3) & 
(({1}^{2},{2}^{1}),2) & 
(({1}^{1},{3}^{1}),1) & 
(({4}^{1}),4) \\ 
\hline
\end{array}
\eeqs
\end{table}
}

Now, we state some new formulas for $a_{5,0}(n)$:
\begin{align}
a_{5,0}(4n) &= a_5(4n), \mylabel{eq:a50form1} \\ 
a_{5,0}(4n+1) &= a_5(4n+1), \mylabel{eq:a50form2} \\ 
a_{5,0}(4n+2) &= 0, \mylabel{eq:a50form3} \\ 
a_{5,0}(4n+3) &= a_5(n). \mylabel{eq:a50form4}  
\end{align}
Formulas (\ref{eq:a50form1})-(\ref{eq:a50form3}) follow from (\ref{eq:elegant1}).
Formula (\ref{eq:a50form4}) is a consequence of the following bijective map,
defined in terms of $n$-vectors by
\beq
\vec{n} \mapsto \vec{n}'
=(2n_1, 1+2n_4, 2n_2, -1+2n_0, 2n_3).
\mylabel{eq:a50trans} 
\eeq
To show that this is a bijection, one may easily verify that
\beq
\abs{\phi_2^{-1}(\vec{n}')} = 4 \abs{\phi_2^{-1}(\vec{n})} + 3,
\mylabel{eq:transprop1} 
\eeq
and show that if $|\pitcb{5}|\equiv3\pmod{4}$ then $\srank(\pitcb{5})\equiv0\pmod{4}$
if and only if
\beq
\phi_2(\pitcb{5}) \equiv (0,1,0,1,0) \pmod{2}.
\mylabel{eq:5csprop}
\eeq

\section{The srank of $t$-cores} \label{sec:sranktcores}

In this section we generalize equation (\ref{eq:elegant1}) to
$t$-cores for general $t$.

\begin{theorem}
\label{theorem4}
Let $t\ge2$, and
\beqs
\phi_2(\pi_{\mbox{$t$-core}}) = \vec{n}=(n_0,n_1,n_2, \dots, n_{t-1}).
\eeqs
Let $a=0$ or $1$.
If $t\equiv 1+2a \pmod{4}$, then
\beq
\srank(\pitc) \equiv \sum_{i=0}^{t-1} (n_i+(1-2a)i+a)^3  \pmod{4};
\mylabel{eq:thm3a}
\eeq
and if $t\equiv 2a \pmod{4}$, then
\beq
\srank(\pitc) \equiv \sum_{i=0}^{t-1} a n_i^2 + (i^2+i)n_i \pmod{4}.
\mylabel{eq:thm3b}
\eeq
\end{theorem}
\begin{proof}
For a partition $\pi=\lamseq$,
\beq
\srank(\pi) \equiv \sum_{j=1}^\nu (\lambda_j^2 + (1-2j)\lambda_j) \pmod{4},
\mylabel{eq:srankcrit}
\eeq
by Proposition 3.1(c) in \cite{Stan}. Now let $t\ge2$, and suppose
$\pi$ is a $t$-core, $\pi=\lamseq$, and
\beqs
\phi_2(\pi) = \vec{n}=(n_0,n_1,n_2, \dots, n_{t-1}).
\eeqs
Suppose $n_i>0$. Then $i$ is exposed in each region $k$, $1\le k \le n_i$,
and this exposed cell is to the right of the main diagonal. Each such
exposed $i$ in region $k$ corresponds to a part of the partition in which
the number of cells to the right of the main diagonal is $t(k-1)+i$.
Thus,
\beq
J = \mbox{size of the Durfee square of $\pi$}
  = \sum_{n_i>0} n_i,
\mylabel{eq:durfee}
\eeq
and
\begin{align}
\srank(\pi)
&\equiv \sum_{j=1}^J (\lambda_j-j)^2 + (\lambda_j - j) + (j - j^2) 
+ \sum_{j=J+1}^\nu \lambda_j^2 + (1-2j) \lambda_j \pmod{4} 
\mylabel{eq:srtc1}\\
&\equiv
\sum_{n_i>0} \sum_{k=1}^{n_i} (t(k-1)+i)^2 + (t(k-1)+i) 
+\sum_{j=1}^J j-j^2 \nonumber\\
&\qquad \qquad
+ \sum_{j=J+1}^\nu \lambda_j^2 + (1-2j) \lambda_j 
\pmod{4} \nonumber \\
&\equiv
\sum_{n_i>0} g(t,n_i,i) + \tfrac{1}{3}(J-J^3)
 + \sum_{j=J+1}^\nu \lambda_j^2 + (1-2j) \lambda_j \pmod{4},
\nonumber
\end{align}
where
\beq
g(t,n,i) = \tfrac{1}{3}t^2n^3 + (ti - \tfrac{1}{2} t(t-1))n^2
+(i^2 - i(t-1) + \tfrac{1}{6}t^2 - \tfrac{1}{2}t)n.
\mylabel{eq:gdef}
\eeq
We note that for $n>0$,
\beq
g(t,n,i) \equiv 0 \pmod{2},
\mylabel{eq:geven}
\eeq
since $g(t,n,i)$ is a sum of even integers. A calculation shows that
\beq
g(t,n,i) + g(t,-n,t-1-i) = 0.
\mylabel{eq:gprop}
\eeq
It follows that $g(t,n,i)$ is a sum of even integers for $n<0$ and \eqn{geven} holds
for all $n$,
so that
\beq
g(t,n,i) \equiv g(t,-n,t-1-i) \pmod{4}.
\mylabel{eq:gmod4}
\eeq
We have
\beq
\phi_2(\pi') = (-n_{t-1},-n_{t-2},\dots,-n_0).
\mylabel{eq:nconj}
\eeq
See \cite[p.3]{GKS}. Let $\pi_1$ be the partition consisting of the first
$J$ parts of $\pi'$. Then by \eqn{nconj}, \eqn{geven} and \eqn{gmod4}
we find that
\begin{align}
\srank(\pi_1)
&\equiv \sum_{n_i<0} g(t,-n_i,t-1-i) + \tfrac{1}{3}(J-J^3) \mylabel{eq:spi1}\\
&\equiv \sum_{n_i<0} g(t,n_i,i)
+ \tfrac{1}{3}(J-J^3) \pmod{4}.\nonumber
\end{align}
Here we have used the fact that $\pi'$ is also a $t$-core.
Now,
\beq
\srank(\pi_1') = - \srank(\pi_1) \equiv \srank(\pi_1) \pmod{4},
\mylabel{eq:spi1c}
\eeq
and
\beq
\srank(\pi_1') \equiv \sum_{j=1}^J J^2 + (1-2j)J
+ \sum_{j=J+1}^\nu \lambda_j^2 + (1-2j) \lambda_j \pmod{4}.
\mylabel{eq:spi1c2}
\eeq
Hence
\beq
\sum_{j=J+1}^\nu \lambda_j^2 + (1-2j) \lambda_j 
\equiv \sum_{n_i<0} g(t,n_i,i) + \tfrac{1}{3}(J-J^3) \pmod{4},
\mylabel{eq:lamsum}
\eeq
since
\beq
\sum_{j=1}^J J^2 + (1-2j)J = 0.
\mylabel{eq:jsum}
\eeq
Hence, by \eqn{srtc1} and \eqn{lamsum}, we have
\begin{align}
\srank(\pi)
&\equiv 
\sum_{n_i>0} g(t,n_i,i)
+\sum_{n_i<0} g(t,n_i,i) + \tfrac{2}{3}(J-J^3) \pmod{4}\mylabel{eq:srtc2}\\
&\equiv 
\sum_{i=0}^{t-1} g(t,n_i,i)\pmod{4},
\nonumber
\end{align}
since $g(t,0,i)=0$ and $\tfrac{2}{3}(J-J^3)\equiv2(J^3-J)\equiv0\pmod{4}$.

We prove \eqn{thm3a} and \eqn{thm3b} 
by finding  simplified forms for $g(t,n,i)\pmod{4}$.
It is clear that the value of 
$g(t,n,i)\pmod{4}$ depends on the residue of $t\pmod{4}$.

\noindent
\underbar{Case 1.} $t\equiv0\pmod{4}$. Then 
\beqs
g(t,n,i)\equiv g(0,n,i)\equiv (i^2+i)n \pmod{4}.
\eeqs
Thus \eqn{thm3b} holds when $a=0$. 

\noindent
\underbar{Case 2.} $t\equiv1\pmod{4}$. Then 
\beqs
g(t,n,i)\equiv g(1,n,i)\equiv -g(1,n,i) \equiv (n+i)^3 - n - i^3 \pmod{4},
\eeqs
so that
\begin{align*}
\srank(\pi)
&\equiv
\sum_{i=0}^{t-1} (n_i+i)^3
-\sum_{i=0}^{t-1} n_i
-\sum_{i=0}^{t-1} i^3 \pmod{4}\\
&\equiv
\sum_{i=0}^{t-1} (n_i+i)^3 - \frac{t^2(t-1)^2}{4} \pmod{4}\\
&\equiv
\sum_{i=0}^{t-1} (n_i+i)^3 \pmod{4},
\end{align*}
since $t\equiv1\pmod{4}$. Here we have also used \eqn{phi2prop1}.
Thus \eqn{thm3a} holds when $a=0$.

\noindent
\underbar{Case 3.} $t\equiv2\pmod{4}$. Then 
\beqs
g(t,n,i)\equiv g(2,n,i) \equiv -g(2,n,i)
\equiv (n^2 + (i^2+i)n) - n + 2in(n+i)
\equiv (n^2 + (i^2+i)n) - n \pmod{4}.  
\eeqs
Thus \eqn{thm3b} holds when $a=1$. 

\noindent
\underbar{Case 4.} $t\equiv3\pmod{4}$. Then 
\beqs
g(t,n,i)\equiv g(3,n,i) \equiv -g(3,n,i)
\equiv (n - i + 1)^3 + n + (i-1)^3 \pmod{4}.
\eeqs
Since   
\beqs
\sum_{i=0}^{t-1} (i-1)^3 = t \frac{(t-3)}{4} ( t(t-3) + 4) 
\equiv 0\pmod{4},
\eeqs
for $t\equiv3\pmod{4}$, we see that \eqn{thm3a} holds when $a=1$.
\end{proof}

\section{The srank and the $t$-quotient} \label{sec:sranktq}       

In this section we prove that
\beq
\mathrm{srank}(\pi) \equiv \mathrm{srank}(\pi_{\mbox{$t$-core}}) 
    + 2a\sum_{i=0}^{t-1} |\pihat_i|\pmod{4},
\mylabel{eq:srtqa}
\eeq
provided $t\equiv 2a\pmod{4}$, and
\begin{align}
\mathrm{srank}(\pi) &\equiv \mathrm{srank}(\pi_{\mbox{$t$-core}}) 
 + 2 \sum_{i=0}^{t-1} (n_i + i + a) \abs{\pihat_i}\mylabel{eq:srtqb}\\ 
&\qquad + \sum_{i=0}^{t-1} \mathrm{srank}(\pihat_i) \pmod{4},
\nonumber
\end{align}
if $t\equiv 1 + 2a \pmod{4}$. Here $a=0$, $1$ and
\beqs
\phi_1(\pi) =(\pitcb{t},(\pihat_0,\dots,\pihat_{t-1})).
\eeqs
We note that when $t=5$ equation \eqn{srtqb} is formula \eqn{elegant2}.
To prove \eqn{srtqa}, \eqn{srtqb} we use \eqn{srankcrit}
which we rewrite as
\beq
\srank(\pi) \equiv \sum_{j=1}^\nu (\lambda_j^2 + (2j-3)\lambda_j) \pmod{4},
\mylabel{eq:srankcritb}
\eeq
where $\pi=\lamseq$.

Next, let $\pi^*$ be a partition obtained from $\pi$ by attaching a single
cell with coordinates $(x,y)$ to the rim of the diagram of $\pi$ as indicated
in Figure \ref{fig2}.

\begin{figure}[tb]
\centerline{\psfig{figure=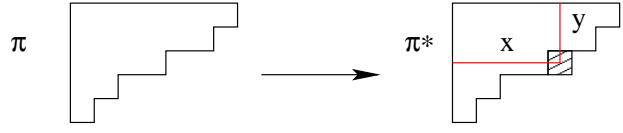}}
\caption{Attaching one cell to the rim}              
\label{fig2}
\end{figure}

It is easy to see that
\beq
\srank(\pi^*)-\srank(\pi)
\equiv(x^2-(x-1)^2) - (y^2-(y-1)^2)\equiv 2(x+y) \pmod{4}.
\mylabel{eq:srpistar}
\eeq
Suppose, we create a new partition $\pi^{**}$ by attaching a border-strip
of length $\ell$ to the diagram of $\pi$, such that the extreme North-East
cell (head) of the strip has coordinates $(x,y)$. Repeated use of \eqn{srpistar}
yields the following formula
\begin{align}
&\srank(\pi^{**})-\srank(\pi)
\equiv
2(x+y) + 2(x+y+1) + \cdots + 2(x+y+\ell-1) \mylabel{eq:srpiss}\\
&\qquad\equiv 2\ell(x+y) + \ell^2 - \ell \pmod{4}.
\nonumber
\end{align}
Here we have used the fact that a border-strip can be added to the diagram one cell
at a time, in such a way that all intermediate diagrams correspond to
partitions. It is straightforward to verify that the right side of \eqn{srpiss}
with $\ell=t\lambda$ becomes
\beq
2a\lambda\pmod{4},\qquad\mbox{if $t\equiv2a\pmod{4}$},\mylabel{eq:srpissa}
\eeq
and
\beq
2a\lambda(x+y+a)+\lambda^2-\lambda\pmod{4},\qquad\mbox{if $t\equiv1+2a\pmod{4}$},
\mylabel{eq:srpissb}
\eeq
with $a=0$, $1$.

Note that \eqn{srpissa} immediately implies \eqn{srtqa}. To prove \eqn{srtqb}
we need to work a little harder. Let's consider a partition $\pitwid_0$,
such that
\beqs
\phi_1(\pitwid_0)
=(\pitcb{t},(\pihat_0,\dots,\pihat_{i-1},(0),\pihat_{i+1},\dots,\pihat_{t-1})).
\eeqs
This partition has the following property. The rim cells with color $i$ are exposed
in all regions $\le n_i$ and are not exposed in all regions $>n_i$ of the
extended $t$-residue diagram of $\pitwid_0$. This means that the word $W_i$ of
$\pitwid_0$ has the form
\beqs
\begin{matrix}
\mbox{Region}: &\cdots\cdots\cdots & n_{i}-1 & n_i & n_{i}+1 & n_{i}+2 & \cdots\cdots\cdots \\
W_i:    &\cdots\cdots\cdots & E       & E   & N       & N       & \cdots\cdots\cdots 
\end{matrix}
\eeqs
Let us attach a border-strip of length $\ell=t\lambda$ to $\pitwid_0$
in such a way that the word $W_i$ becomes
\beqs
\begin{matrix}
\mbox{Region}: &\cdots & n_{i}-\lambda_1 & n_i+1-\lambda_1 & \cdots & \cdots & n_{i}+1 
& n_{i}+2 & n_i + 3 & \cdots \\
W_i:    &\cdots & E               & N               & E      & \cdots & E 
& N       & N       & \cdots 
\end{matrix}
\eeqs
This way we create a partition $\pitwid_1$ such that
\beqs
\phi_1(\pitwid_1)
=(\pitcb{t},(\pihat_0,\dots,\pihat_{i-1},(\lambda_1),\pihat_{i+1},\dots,\pihat_{t-1})).
\eeqs
It is straightforward to verify that for $t$ odd the coordinates $(x,y)$ of the
border-strip head satisfy
\beq
x + y \equiv n_i + i \pmod{2}.
\mylabel{eq:bsh}
\eeq
Hence using \eqn{srpissb}, \eqn{bsh} we find that for $t\equiv1+2a\pmod{4}$
\beq
\srank(\pitwid_1)-\srank(\pitwid_0)
\equiv 2\lambda_1(n_i + i + a) + \lambda_1^2 - \lambda_1\pmod{4},
\mylabel{eq:srq1}
\eeq
where $a=0$, $1$.

Next, we add to the diagram of $\pitwid_1$ a new border-strip of length $t\lambda_2$
with $\lambda_2\le\lambda_1$, such that $W_i$ becomes

\beqs
\begin{matrix}
\mbox{Region}: &\cdots \cdots &n_{i}+1-\lambda_1 &\cdots \cdots \cdots     &n_{i}+2 -\lambda_2 
&\cdots \cdots \cdots &n_i + 3 &\cdots\cdots\\ 
W_i:    &\cdots \quad E& N                & E \quad\cdots \quad E   &N    
&E \quad\cdots\quad E &N        &N\quad\cdots
\end{matrix}
\eeqs
The new partition $\pitwid_2$ satisfies
\beqs
\phi_1(\pitwid_2)
=(\pitcb{t},(\pihat_0,\dots,\pihat_{i-1},(\lambda_1,\lambda_2),\pihat_{i+1},\dots,\pihat_{t-1})).
\eeqs
Replacing $\lambda_1$ by $\lambda_2$ and $n_i$ by $n_i+1$ and repeating the
argument that led us to \eqn{srq1} we obtain
\beq
\srank(\pitwid_2)-\srank(\pitwid_1)
\equiv 2\lambda_2(n_i + 1 + i + a) + \lambda_2^2 - \lambda_2\pmod{4}.
\mylabel{eq:srq2}
\eeq
Let $\pitwid_\nu$ denote the partition obtained from $\pitwid_0$, such that
\beqs
\phi_1(\pitwid_\nu)
=(\pitcb{t},(\pihat_0,\dots,\pihat_{i},\dots,\pihat_{t-1})),
\eeqs
where $\pihat_i=\lamseq$. Proceeding as above and using \eqn{srankcritb}
we find that for $t\equiv1 +2a\pmod{4}$
\begin{align}
\srank(\pitwid_\nu)-\srank(\pitwid_0) &
\equiv 2 \sum_{j=1}^\nu \lambda_j (n_i + i + a + j-1) 
+  \sum_{j=1}^\nu (\lambda_j^2 - \lambda_j) \mylabel{eq:srq3}\\
& \equiv 2(n_i + i + a ) \sum_{j=1}^\nu \lambda_j  
+  \sum_{j=1}^\nu (\lambda_j^2  + (2j-3)\lambda_j) \nonumber \\
& \equiv 2(n_i + i + a ) |\pihat_i| + \srank(\pihat_i) \pmod{4}, 
\nonumber
\end{align}
with $a=0$, $1$. Formula \eqn{srtqb} follows easily from \eqn{srq3}.

\section{Generalization of Andrews' refinement and new partition congruences modulo $5$} 
\label{sec:bgrank}        
In Section \ref{subsec:2quotrank} we gave a new combinatorial interpretation
of Andrews' result \eqn{andrefine} in terms of the $2$-quotient of a partition.
Further study of this development led us to a generalization of \eqn{andrefine},
which we now describe. We define the new partition statistic
\beq
\mbgr(\pi) = \sum_{j=1}^\nu (-1)^{j+1} \parity(\lambda_j),
\mylabel{eq:bgrdef}
\eeq
where $\pi=\lamseq$ and for an integer $m$, $\parity(m)$ denotes the
parity of $m$;\ i.e.\ $\parity(m)=1$ if $m$ is odd and $0$, otherwise.
If $\phi_1(\pi)=(\pitcb{2},(\pihat_0,\pihat_1))$ and $\phi_2(\pitcb{2})=(n_0,-n_0)$
then it is easy to verify that
\beq
\mbgr(\pi) = n_0 = r_0 - r_1.
\mylabel{eq:bgralt}
\eeq
Here $r_i$ with $i=0$, $1$ denotes the number of cells colored $i$ in the
$2$-residue diagram of $\pi$.

Next, we recall that
\beq
\srank(\pi) \equiv |\pi| - |\pitcb{2}|
\equiv |\pi| - n_0(2n_0 -1) \pmod{4}.
\mylabel{eq:srn0}
\eeq
Here we have used \eqn{sr2core} and \eqn{phi2prop2}. Thus, if $|\pi|$
is given, then $\srank(\pi)\pmod{4}$ is completely determined by
$\mbgr(\pi)$. Clearly, the converse is not true.

Let $\pbar_{j}(m,n)$ denote the number of partitions of $n$ with 
$\mbgr=j$ and $\mtqr=m$. Then
\begin{align}
f_j(x,q) &=
\sum_{\substack{n\ge0,\\ m}} \pbar_j(m,n) x^m q^n \mylabel{eq:fjdef}\\
&= q^{(2j-1)j} \sum_{\pihat_0,\pihat_1} q^{2(|\pihat_0|+|\pihat_1|)}
x^{\nu(\pihat_0)-\nu(\pihat_1)}\nonumber\\
&= 
\frac{q^{(2j-1)j}}{(q^2x,q^2/x;q^2)_\infty}.
\nonumber
\end{align}
Proceeding as in Section \ref{sec:stcrank}, we find that for $\xi^5=1$,
$\xi\ne1$
\beq
f_j(\xi,q) = 
\frac{1}{1-\xi^2} \frac{q^{(2j-1)j}}{(q^{10};q^{10})_\infty}
\sum_{n\ge0} (-1)^n q^{n^2+n} \xi^{-2n}(1-\xi^{4n+2}).
\mylabel{eq:fjxi}
\eeq
We have
\beq
f_j(\xi,q) =
\sum_{\substack{n\ge0\\ 0\le m \le4}} \Pbar_j(m,5,n) \xi^m q^n,
\mylabel{eq:fjxi2}
\eeq
where $\Pbar_j(m,5,n)$ denotes the number of partitions of $n$ with
$\mbgr=j$ and $\mtqr\equiv m\pmod{5}$.

Next, we observe that
\beq
n^2 + n \equiv
\begin{cases}
0\pmod{5}, &\mbox{if $n\equiv0,4\pmod{5}$},\\
2\pmod{5}, &\mbox{if $n\equiv1,3\pmod{5}$},\\
1\pmod{5}, &\mbox{if $n\equiv2\pmod{5}$},
\end{cases}
\mylabel{eq:trim5}
\eeq
\beq
(2j-1)j \equiv
\begin{cases}
0\pmod{5}, &\mbox{if $j\equiv0,3\pmod{5}$},\\
1\pmod{5}, &\mbox{if $j\equiv1,2\pmod{5}$},\\
3\pmod{5}, &\mbox{if $j\equiv4\pmod{5}$},
\end{cases}
\mylabel{eq:2cm5}
\eeq
and $(1-\xi^{4n+2})=0$ if $n\equiv2\pmod{5}$. This means that
\begin{align}
\sum_{m=0}^4 \Pbar_j(m,5,5n) \xi^m &=0, \quad\mbox{if $j\equiv1,2\pmod{5}$},
\mylabel{eq:Pjrel1}\\
\sum_{m=0}^4 \Pbar_j(m,5,5n+1) \xi^m &=0, \quad\mbox{if $j\not\equiv1,2\pmod{5}$},
\mylabel{eq:Pjrel2}\\
\sum_{m=0}^4 \Pbar_j(m,5,5n+2) \xi^m &=0, \quad\mbox{if $j\not\equiv0,3\pmod{5}$},
\mylabel{eq:Pjrel3}\\
\sum_{m=0}^4 \Pbar_j(m,5,5n+3) \xi^m &=0, \quad\mbox{if $j\equiv0,3\pmod{5}$},
\mylabel{eq:Pjrel4}\\
\sum_{m=0}^4 \Pbar_j(m,5,5n+4) \xi^m &=0, \quad\mbox{for all $j$}.
\mylabel{eq:Pjrel5}
\end{align}
Hence we have the following
\begin{theorem}
\label{theorem5}
For $m=0$, $1$, $2$, $3$, $4$
\begin{align}
\Pbar_j(m,5,5n) &= \frac{1}{5} \pbar_j(5n), \quad\mbox{if $j\equiv1,2\pmod{5}$},
\mylabel{eq:thm5a}\\
\Pbar_j(m,5,5n+1) &= \frac{1}{5} \pbar_j(5n+1), \quad\mbox{if $j\not\equiv1,2\pmod{5}$},
\mylabel{eq:thm5b}\\
\Pbar_j(m,5,5n+2) &= \frac{1}{5} \pbar_j(5n+2), \quad\mbox{if $j\not\equiv0,3\pmod{5}$},
\mylabel{eq:thm5c}\\
\Pbar_j(m,5,5n+3) &= \frac{1}{5} \pbar_j(5n+3), \quad\mbox{if $j\equiv0,3\pmod{5}$},
\mylabel{eq:thm5d}\\
\noalign{\text{and}}
\Pbar_j(m,5,5n+4) &= \frac{1}{5} \pbar_j(5n+4), \quad\mbox{for all $j$}.
\mylabel{eq:thm5e}
\end{align}
\end{theorem}
\begin{corollary}
\label{cor5}
\begin{align}
\pbar_j(5n) &\equiv0\pmod{5} \quad\mbox{if $j\equiv1,2\pmod{5}$},
\mylabel{eq:cor5a}\\
\pbar_j(5n+1) &\equiv0\pmod{5} \quad\mbox{if $j\not\equiv1,2\pmod{5}$},
\mylabel{eq:cor5b}\\
\pbar_j(5n+2) &\equiv0\pmod{5} \quad\mbox{if $j\not\equiv0,3\pmod{5}$},
\mylabel{eq:cor5c}\\
\pbar_j(5n+3) &\equiv0\pmod{5} \quad\mbox{if $j\equiv0,3\pmod{5}$},
\mylabel{eq:cor5d}\\
\noalign{\text{and}}
\pbar_j(5n+4) &\equiv0\pmod{5} \quad\mbox{for all $j$}.
\mylabel{eq:cor5e}
\end{align}
\end{corollary}
Recalling the comment after \eqn{srn0} we see that \eqn{cor5e} gives
an extension of Andrews' result \eqn{andrefine}.
However, congruences \eqn{cor5a}--\eqn{cor5d} appear to be new.
It is possible to modify the construction in Section \ref{sec:refinement}
in order to provide a direct combinatorial proof of \eqn{cor5e}.
It is likely that a combinatorial proof of \eqn{cor5a}--\eqn{cor5d}
would require significant new insights. 
To make this point plausible we note that the $5$-core analog of
\eqn{cor5a}--\eqn{cor5d} does not hold. In other words, it is not true
that for $r=0$, $1$, $2$, $3$      
\beq
\abar_{5,j}(5n+r) \equiv 0 \pmod{5}.
\mylabel{eq:ab5jr}
\eeq
On the other hand, it can be shown that
\beq
\abar_{5,j}(5n+4) \equiv 0 \pmod{5}.
\mylabel{eq:ab5j4}
\eeq
Here $\abar_{5,j}(n)$ denotes the number of $5$-core partitions
of $n$ with $\mbgr=j$.

Unfortunately, neither the srank nor the \bgr\ can be 
employed to refine Ramanujan's congruences mod $7$ or $11$. And so we would
like to pose the 

\noindent
\textbf{Problem}. Is there an analogue of the \bgr, which gives a        
refinement of \eqn{ram7} and \eqn{ram11}?

{\it Acknowledgements}.
We would like to thank Krishna Alladi, George Andrews, and Alain Lascoux for 
their interest and comments.

\bibliographystyle{amsplain}

\begin{thebibliography}{10}
%
\bibitem{Andrews1}
G. E. Andrews,
\textit{On a partition function of Richard Stanley}, to appear in the
Electronic Journal of Combinatorics volume in honor of Richard
Stanley.
%
\bibitem{AG}
G. E. Andrews and F. G. Garvan,
\textit{Dyson's crank of a partition},
{Bull. Amer. Math. Soc. (N.S.)}
\textbf{18} (1988), 167--171.
%
\bibitem{ASD}
A.~O.~L.~Atkin and P.~Swinnerton-Dyer,
\textit{Some properties of partitions}, 
Proc. London Math. Soc.
\textbf{4} (1954), 84--106.
%
\bibitem{Boulet}
C. E. Boulet,
\textit{A four-parameter partition identity},
preprint.           
%
\bibitem{Dyson}
F. J. Dyson,
\textit{Some guesses in the theory of partitions},
Eureka (Cambridge)
\textbf{8} (1944), 10--15.
%
\bibitem{G1}
F. G. Garvan,
\textit{New combinatorial interpretations of Ramanujan's partition
congruences mod $5,7$ and $11$},
Trans. Amer. Math. Soc.
\textbf{305} (1988), 47--77.
%
\bibitem{G2}
F. G. Garvan,
\textit{The crank of partitions mod $8,\;9$ and $10$},
{Trans. Amer. Math. Soc.}
\textbf{322} (1990), 79--94.
%
\bibitem{G3}
F. G. Garvan,
\textit{More cranks and $t$-cores},
Bull. Austral. Math. Soc.
\textbf{63} (2001), 379--391.
%
\bibitem{GKS}
F. Garvan, D. Kim and D. Stanton,
\textit{Cranks and $t$-cores},
{Invent. Math.}
\textbf{101} (1990), 1--17.
%
\bibitem{HGB}
M. Hirschhorn, F. Garvan and J.Borwein,
\textit{Cubic analogues of the Jacobian theta function $\theta(z,q)$},
{Canad. J. Math.}
\textbf{45} (1993), 673--694.  
%
\bibitem{JK}
G. James and A. Kerber,
\textit{The Representation Theory of the Symmetric Group},
Addison-Wesley, Reading, MA, 1981.
%
\bibitem{L}
D. E. Littlewood, 
\textit{Modular representations of symmetric groups},
{Proc. Roy. Soc. London. Ser. A.}
\textbf{209} (1951), 333--353.
%
\bibitem{Sills}
A. V. Sills,
\textit{A combinatorial proof of a partition identity of
Andrews and Stanley},
preprint.
%
\bibitem{Stan}
R. P. Stanley,
\textit{Some remarks on sign-balanced and maj-balanced posets},
preprint.
%
\bibitem{Y}
A. J. Yee,
\textit{On partition functions of Andrews and Stanley},
preprint.
%
\end{thebibliography}

\end{document}